\documentclass[journal,twocolumn,10pt,twoside]{IEEEtran}

\usepackage{amsmath,amsfonts, amssymb,color}
\usepackage{epstopdf}
\usepackage{nomencl}
\usepackage{enumerate}
\usepackage{times} 
\usepackage{latexsym}
\usepackage{xspace}
\usepackage{amssymb,amsmath,amsfonts,color}
\usepackage{cite}
\usepackage{caption}
\usepackage{tikz,subcaption,cite}
\usetikzlibrary{decorations.pathmorphing}
\usepackage{soul}

\newtheorem{problem}{Problem}
\newtheorem{theorem}{Theorem}

\newtheorem{lemma}{Lemma}
\newtheorem{remark}{Remark}

\newtheorem{definition}{Definition}
\newtheorem{example}{Example}

\renewcommand{\baselinestretch}{1.05} 

\definecolor{cyan}{rgb}{0.0, 1.0, 1.0}

\begin{document}
\title{Optimal Active Fault Detection in\\ Inverter-Based Grids}


\author{Mohammad Pirani, Mehdi Hosseinzadeh, \IEEEmembership{Member, IEEE}, Joshua A. Taylor, and Bruno Sinopoli, \IEEEmembership{Fellow, IEEE}
\thanks{M. Pirani and J. A. Taylor are with the Department of Electrical and Computer Engineering, University of Toronto. E-mail: {\tt mohammad.pirani@utoronto.ca; josh.taylor@utoronto.ca}.} 
\thanks{M. Hosseinzadeh is with the School of Mechanical and Materials Engineering, Washington State University, Pullman, WA 99164, USA. E-mail: {\tt mehdi.hosseinzadeh@wsu.edu}.}
\thanks{B. Sinopoli is with the Department of Electrical and Systems Engineering, Washington University in St Louis, St. Louis, MO 63130, USA. E-mail: {\tt bsinopoli@wustl.edu}.}%
}

\maketitle


\begin{abstract}
Ground faults in converter-based grids can be difficult to detect because, unlike in grids with synchronous machines, they often do not result in large currents. One recent strategy is for each converter to inject a perturbation that makes faults easier to distinguish from normal operation. In this paper, we construct optimal perturbation sequences for use with the Multiple Model Kalman Filter. The perturbations maximize the difference between faulty and fault-free operation while respecting limits on performance degradation.  Simulations show that the optimal input sequence increases the confidence of fault detection while decreasing detection time. It is shown that there is a tradeoff between detection and degradation of the control performance, and that the method is robust  to parameter variations. 
\end{abstract}




\section{Introduction}
\subsection{Motivation}
Faults are relatively easy to detect in conventional power grids because of the large currents they draw from synchronous machines. In inverter-based grids, e.g., microgrids, the converter controllers prevent fault currents from rising substantially above nominal \cite{hooshyar, surveymicro}. As a result, traditional detection-schemes might fail in inverter-based grids.

One promising approach is for each converter to inject a small perturbation that makes faults easier to detect without degrading quality. This was done in \cite{karimi} to detect unintentional islanding, and more recently in \cite{khaled} to detect faults. We build on the latter by optimizing the perturbation sequence. This maximizes the probability that a multiple model Kalman filter (MMKF) \cite{vhassani} will be able to distinguish between faulty and fault-free operation.

\subsection{Literature Review}
A fault in a power system is an unintentional electrical path to the ground or between phases.  The standard remedy is to  de-energize them by opening protection relays \cite{Phadke}. Detection of faults as fast as possible is vital to protect the grid. Several passive, e.g., signal processing techniques \cite{dong, Borghetti, Lorenc}, and active methods, e.g., zero and negative sequence power injection  \cite{lin, karimi}, have been proposed in the literature.

The current controllers in inverters prevent faults from drawing large currents in inverter-based grids, which makes faults harder to detect.  To improve detectability, \cite{khaled} proposed injecting a harmonic signal to increase the difference between normal and faulty operation.  We build on this by injecting an optimized signal instead of a harmonic, and by using an MMKF to detect faults.

The MMKF determines which mode of operation a system is in. In this case, it determines if the system is faulty or not. The MMKF  is inspired from multiple-model adaptive estimation \cite{athans, li, hasani}.  More recently the MMKF was used to detect faults and attacks \cite{Mehdi}. We apply this to faults in inverter-based grids.

\subsection{Contributions}
In this paper, we use the MMKF to design an active fault detection scheme for inverter-based grids. We evaluate the scheme in simulation, and find that the optimal perturbation sequence dramatically enhances detection, and outperforms harmonic-based schemes.

We remark that here we only consider three phase, balanced faults and a single converter. Most real faults, however, are unbalanced, and most grids have multiple converters and synchronous machines. Extending our approach to more realistic settings is a topic of future work.

The structure of the paper is as follows. In Section \ref{sec:nc092}, model the system. In Section \ref{sec:cn09h1}, we present the MMKF, show how to solve for the optimal input sequence, and discuss their application to fault detection. Section \ref{sec:cm09n72} has the simulation results and Section \ref{sec:conclusion} concludes the paper.

\section{Modeling}
\label{sec:nc092}
We consider a single inverter-based distributed generator (DG). We assume that the system is balanced and only consider balanced, three phase-to-ground faults. We model the system in the $dq$ frame.

\subsection{Inverter Modeling}
 The control objective is to keep the voltage of the DG at some nominal value. This is commonly achieved with an inner controller that keeps the current in a safe range, and an outer controller that tracks the voltage setpoint and sets the reference for the inner controller \cite{Jovcic}.

\begin{figure}[!t]
\centering
\includegraphics[scale=0.58]{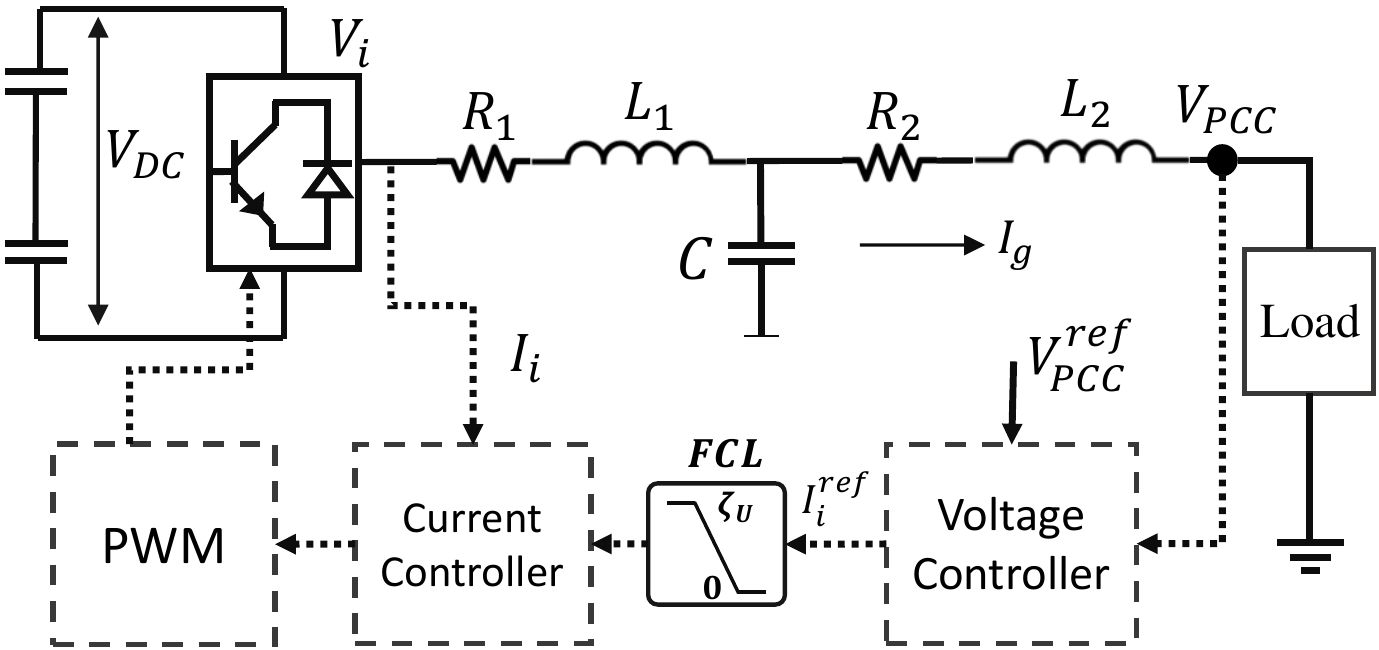}
\caption{ Inverter and control loops.}
\label{fig:exwerae2ewdmple}
\end{figure}
The system is shown in  Fig.~\ref{fig:exwerae2ewdmple}. $V_i$ is the inverter  output voltage, $V_{\rm DC}$  the  DC voltage, and $I_g$  the grid-side line current.  $R_1, L_1$ and $R_2, L_2$ are converter and the transformer's impedances, respectively. $C$ and $I_C$ are the filter capacitance and its current. The voltage controller drives $V_{\rm PCC}$ to a reference value, and sends a reference signal to the current controller.  The inverter current must track this reference value. The fault current limiter (FCL) is a saturation block that protect the converter from large currents \cite{yazdani}. 

\subsection{Voltage and Current Control Loops}
Let $e_V(t)=V_{\rm PCC}^{\rm ref}-V_{\rm PCC}$. The voltage at the PCC is managed by the PI controller
\begin{equation}\label{eqn:control1}
u_V(t)=k_{p}^Ve_V(t)+k_{i}^V\int_0^te_V(\tau)d\tau,
\end{equation} 
where  $k_{p}^V$ and $k_{i}^V$ are  gains. The voltage controller determines the reference inverter  current. The current controller is  designed starting from the inverter's AC voltage \cite{Jovcic}, which in the $dq$ frame is given by  
\begin{equation}\label{eqn:hvs8hobiuaf}
    V_{i}^d=\frac{1}{2}V_{\rm DC}u_{I}^d,  \quad  V_{i}^q=\frac{1}{2}V_{\rm DC}u_{I}^q,
\end{equation}
where $u_{I}^d$ and $u_{I}^q$ are duty cycles. 
Denoting  the capacitor voltage by $V_C$, the dynamics of the inverter line current are
\begin{equation}\label{eqn:cn0898586a}
V_C^j=V_i^j-R_1I_i^j-L_1\frac{dI_i^j}{dt}, \quad \quad j\in\{a,b,c\}.
\end{equation}
Transforming \eqref{eqn:cn0898586a} into the $dq$ frame and using \eqref{eqn:hvs8hobiuaf} gives
\begin{align}\label{eqn:c568ssbg86a}
V_C^d&=\frac{1}{2}V_{\rm DC}u_{I}^d-R_1I_i^d+\omega L_1I_i^q-L_1\frac{dI_i^d}{dt}, \nonumber\\
V_C^q&=\frac{1}{2}V_{\rm DC}u_{I}^q-R_1I_i^q-\omega L_1I_i^d-L_1\frac{dI_i^q}{dt}.
\end{align}
To decouple \eqref{eqn:c568ssbg86a}, define the new control signals
\begin{align}
u_{I}^d=2\frac{\hat{u}_{I}^d-L_1\omega I_i^q}{V_{\rm DC}},\quad
u_{I}^q=2\frac{\hat{u}_{I}^q+L_1\omega I_i^d}{V_{\rm DC}}.
\end{align}
The dynamics are now decoupled and take the form 
\begin{align}\label{eqn:c56168sd4gaga}
V_C^d&=\hat{u}_{I}^d-R_1I_i^d-L_1\frac{dI_i^d}{dt}, \nonumber\\
V_C^q&=\hat{u}_{I}^q-R_1I_i^q-L_1\frac{dI_i^q}{dt}.
\end{align}
The current controller is given by
$$\hat{u}_I^j(t)=k_{p}^Ie_I^j(t)+k_{i}^I\int_0^te_I^j(\tau)d\tau, \quad \quad j\in \{d,q\}$$ 
where $e_I^d(t)=I_{i}^{\rm ref,d}-I_{i}^d$ and $e_I^q(t)=I_{i}^{\rm ref,q}-I_{i}^q$.

\subsection{Fault Modeling}
\begin{figure}[t!]
\centering
\includegraphics[scale=0.54]{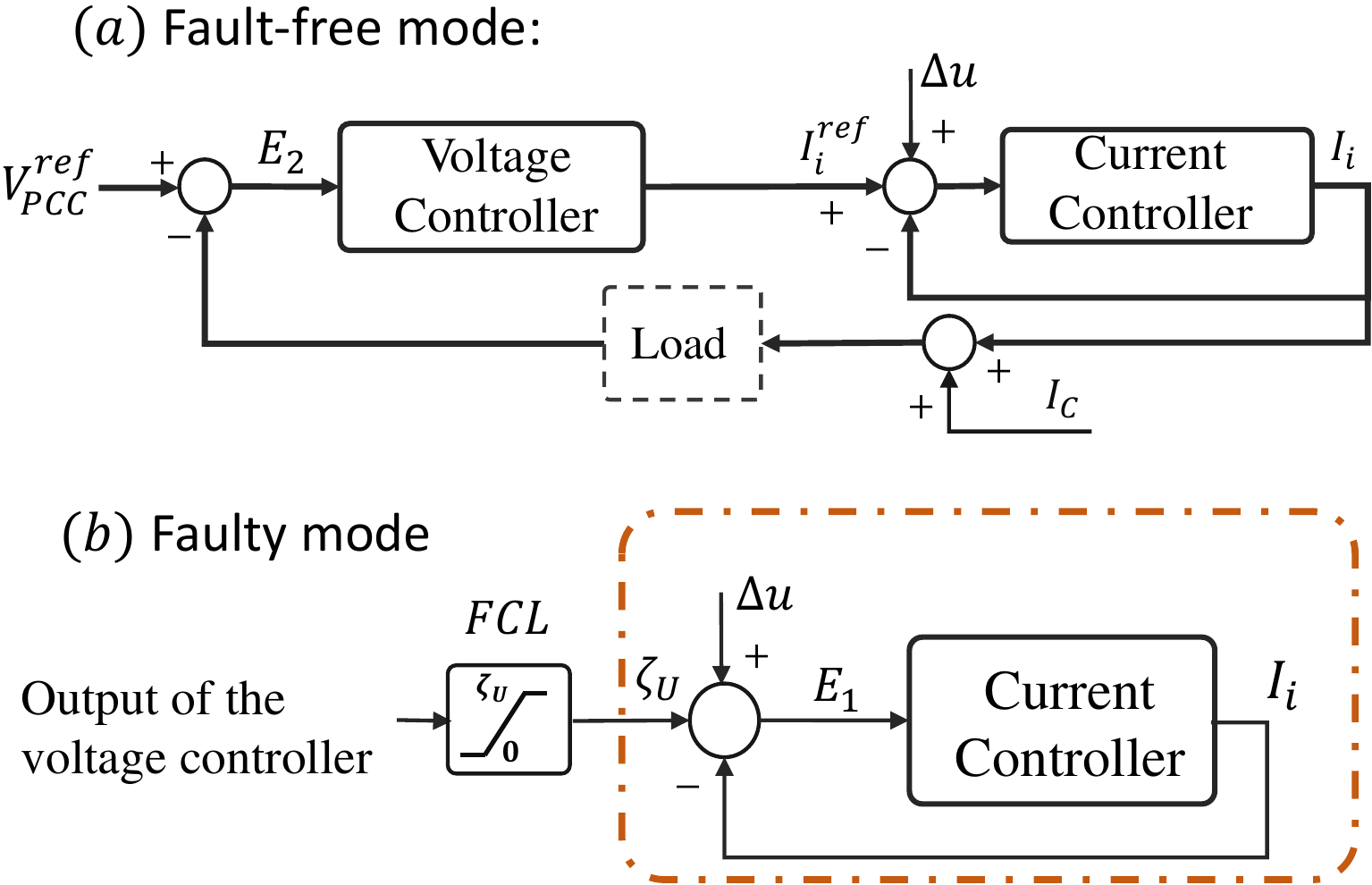}
\caption{(a) Fault-free mode, in which the FCL is not active.   (b) Faulty mode in which the FCL is active.}
\label{fig:exwsghrle}
\end{figure}

We focus on three phase symmetrical ground faults. The controller during fault-free operation is shown in Fig.~\ref{fig:exwsghrle} (a).
In  this case, the FCL is not active. 
The control loop during a ground fault is shown in Fig.~\ref{fig:exwsghrle} (b). In this case the current increases, and the upper limit of the FCL, $\zeta_U$, becomes the reference for the current controller.  The capacitor $C$ is designed to compensate the current deficiency in the converter side during faults. 
\begin{remark}
The switching in the FCL is much faster than the dynamics of the controllers. Due to this time scale separation, we do not account for the transient aspects of saturation in the fault detection algorithm. 
\end{remark}

\subsection{State Space Representation}

During fault-free operation, the system has the following model.
\begin{align}\label{eqn:cn09j2svs}
\dot{x}&=I_{2}\otimes A_{h}x +I_{2}\otimes B_{h} u
\nonumber\\
y&=\begin{bmatrix}
0 &1 & 0 & 0 & 0 & 0\\
0 &0 & 0 & 0 & 1 & 0
\end{bmatrix}x
\end{align}
where $I_2$ is the $2\times2$ identity matrix and $\otimes$ is the Kronecker product and
\begin{align}\label{eqn:cn09j2}
A_{h}&= \begin{bmatrix}
0 & -k_{i}^I-k_{i}^Ik_{p}^V & k_{i}^I\\ \frac{V_{\rm DC}}{2L_1} & -\left(\frac{k_{p}^IV_{\rm DC}}{2L_1}+\frac{R_1}{L_1}+\frac{k_{p}^Ik_{p}^VV_{\rm DC}R}{2L_1}\right) & \frac{k_{p}^IV_{\rm DC}}{2L_1} \\
0 & -k_{i}^VR & 0
\end{bmatrix}, \nonumber\\ B_{h}&=\begin{bmatrix}
k_{i}^Ik_{p}^V & k_{i}^I\\ \frac{k_{p}^Ik_{i}^IV_{\rm DC}}{2L_1} & \frac{k_{p}^IV_{\rm DC}}{2L_1} \\ k_{i}^V & 0
\end{bmatrix}.
\end{align}
Here, $x=[x_1,x_2,x_3,x_4,x_5,x_6]^T$ where $x_2=I_i^d$, $x_5=I_i^q$, and $x_1, x_3,x_4,x_6$ are intermediate states in the integral paths of PI controllers, i.e., $\dot{x}_1=k_{i}^IE_{1}^d$, $\dot{x}_3=k_{i}^IE_{1}^q$, $\dot{x}_4=k_{i}^VE_{2}^d$, $\dot{x}_6=k_{i}^VE_{2}^q$,
where $E_1$ and $E_2$ are the errors shown in Fig.~\ref{fig:exwerae2ewdmple} (a) and (b). 
 The input is $u=\begin{bmatrix}
V_{\rm PCC}^{\rm ref,d} & \Delta u^d & V_{\rm PCC}^{\rm ref,q} & \Delta u^q
\end{bmatrix}^T$ where $\Delta u$ is a perturbation (more precisely, $\Delta u$ is a variation in the reference signal of the current controller) which we describe later. 
The output is the inverter current, $I_i$, i.e., states $x_2$ and $x_5$. In the above dynamics, we have assumed a resistive load with equivalent resistance $R$.

During a fault, the model is given by
\begin{align}
\dot{x}&=I_{2}\otimes A_{f}x +I_{2}\otimes B_{f} u\nonumber\\
y&=\begin{bmatrix}
0 &1&0&0\\0 &0&0&1
\end{bmatrix}x
\end{align}
where
\begin{align}\label{eqn:cn09j2}
A_{f}=\begin{bmatrix}
0 & -k_{i}^I\\ \frac{V_{\rm DC}}{2L_1} & -\left(\frac{k_{p}^IV_{\rm DC}}{2L_1}+\frac{R_1}{L_1}\right)
\end{bmatrix}, B_{f}=\begin{bmatrix}
k_{i}^I & k_{i}^I\\ \frac{k_{p}^IV_{\rm DC}}{2L_1} & \frac{k_{p}^IV_{\rm DC}}{2L_1} 
\end{bmatrix}.
\end{align}
Here, $x=[x_1,x_2,x_3,x_4]^T$ where $x_2=I_i^d$, $x_4=I_i^q$, and $x_1, x_3$ are intermediate states with $\dot{x}_1=k_{i}^IE_{1}^d$ and $\dot{x}_3=k_{i}^IE_{1}^q$. The input in this case is $u=\begin{bmatrix}
\zeta_U^d & \Delta u^d & \zeta_U^q & \Delta u^q
\end{bmatrix}^T$ and, similar to the fault-free case, the output is the inverter current. In the rest of the paper, $\Delta u(t):=[\Delta u^d(t)~\Delta u^q(t)]^\top$.

\section{Multiple Model Kalman Filter}
\label{sec:cn09h1}

We first review the MMKF. Let ${h}$ and $f$ be the indices of fault-free  and faulty modes, respectively. We discretize and combine the fault-free and faulty modes as 
\begin{align}
x_{k+1|i}&=A_{i}x_{k|i}+B_{i}u_{k}+w_{k},   \nonumber\\
 y_{k|i}&=C_{i}x_{k|i}+v_{k},
 \label{eqn:cna082}
\end{align}
where $i\in\{h,f\}$, and $x_{k|h}\in\mathbb{R}^{n_{h}}$ and  $x_{k|f}\in\mathbb{R}^{n_{f}}$ are respectively the state vectors of fault-free and faulty systems at time $k$, with $n_h>0$ and $n_f>0$ as the dimension of the state vectors. Here, $v_{k}$ and $w_{k}$  are measurement and process noises with covariances $\Sigma_v$ and $\Sigma_w$, respectively. Both modes start from the  initial condition  $x_{0|i}=x_0$ with covariance $\Sigma_0$.

\begin{figure}[t!]
\centering
\includegraphics[scale=0.65]{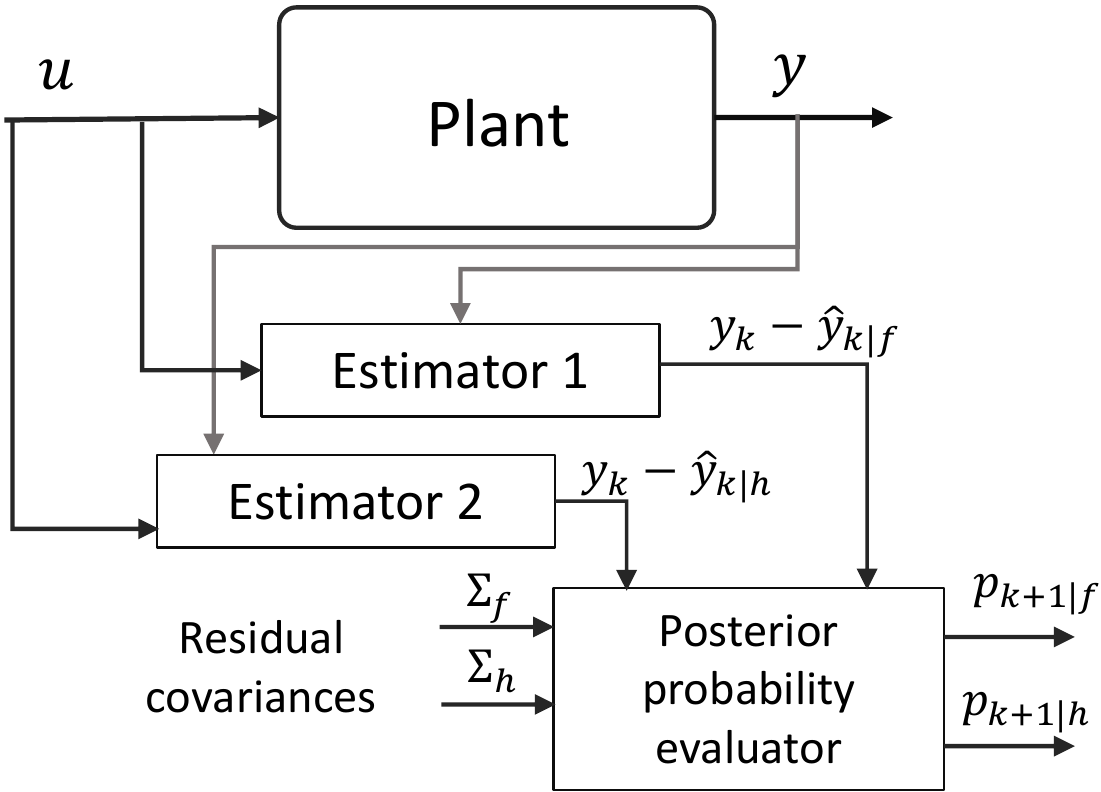}
\caption{Multiple-Model Kalman Filter.}
\label{fig:exwere9emple}
\end{figure}

The MMKF is  shown in Fig.~\ref{fig:exwere9emple}. The procedure starts from an initial guess of the probabilities of the system being faulty or fault-free, denoted by $p_{0|h}$ and $p_{0|f}$. Then, we estimate the output of the system for each mode and update the probability based on the residual signals, i.e., the difference between the measured output and the estimated output. More formally, using measurements at each time step $k$, the Kalman filter for each mode $i\in \{h,f\}$ is given by
\begin{align}
\hat{x}_{k+1|i}&=A_{i}\hat{x}_{k|i}+B_{i}u_k+H_{i}\left( y_k-C_{i}\hat{x}_{k|i} \right)\nonumber\\
H_{i}&=\Sigma_{i}C_{i}^T\left[C_{i}\Sigma_{i}C_{i}^T+\Sigma_v  \right]^{-1}
\end{align}
where $\hat{x}_{k|h}\in \mathbb{R}^{n_{h}}$ and  $\hat{x}_{k|f}\in \mathbb{R}^{n_{f}}$ are the estimated states at time $k$, respectively. $\Sigma_{i}$ is the posterior error covariance matrix in steady state,  which is the  solution of the discrete Riccati equation, 
\begin{align}
0&=-\Sigma_{i}+A_{i}\Sigma_{i}A_{i}^T\nonumber\\
&+ \Sigma_w- A_{i}^T\Sigma_{i}C_{i}^T[C_{i}\Sigma_{i}C_{i}^T+\Sigma_v]^{-1}C_{i}\Sigma_{i}A_{i}.
\label{eqn:n0j9ahc}
\end{align}

After estimating the states of each mode, 
we use the measurement $y_k$ to generate the  residual functions $\alpha_{k|f}=\|y_k-\hat{y}_{k|f}\|$ and $\alpha_{k|h}=\|y_k-\hat{y}_{k|h}\|$.  These residuals are used to update the posterior probability of each mode
\begin{equation}\label{eqn:posterior}
p_{k+1|i}=\frac{\beta_ie^{-\alpha_{k|i}}}{p_{k|f}\beta_fe^{-\alpha_{k|f}}+p_{k|h}\beta_he^{-\alpha_{k|h}}}  p_{k|i},
\end{equation}
where $\beta_i=\frac{1}{\sqrt{|\Sigma_{y|i}|}}$ and $\Sigma_{y|i}$ is the  covariance of the residual $y_k-\hat{y}_{k|i}$. The residual is  updated as
\begin{align}
\Sigma_{y|i}=C_{i}\Sigma_{i}C_{i}^T+\Sigma_{v}
\end{align}
with $\Sigma_{i}$ as in \eqref{eqn:n0j9ahc}.

The detected mode is the one for which the corresponding posterior probability is maximized, i.e., 
\begin{equation}
    p^*=\arg\max_{i\in \{h,f\}}p\left(i|y_{0:N},\Delta u_{0:N-1}\right)
\end{equation}

\begin{figure}[t!]
\centering
\includegraphics[scale=0.58]{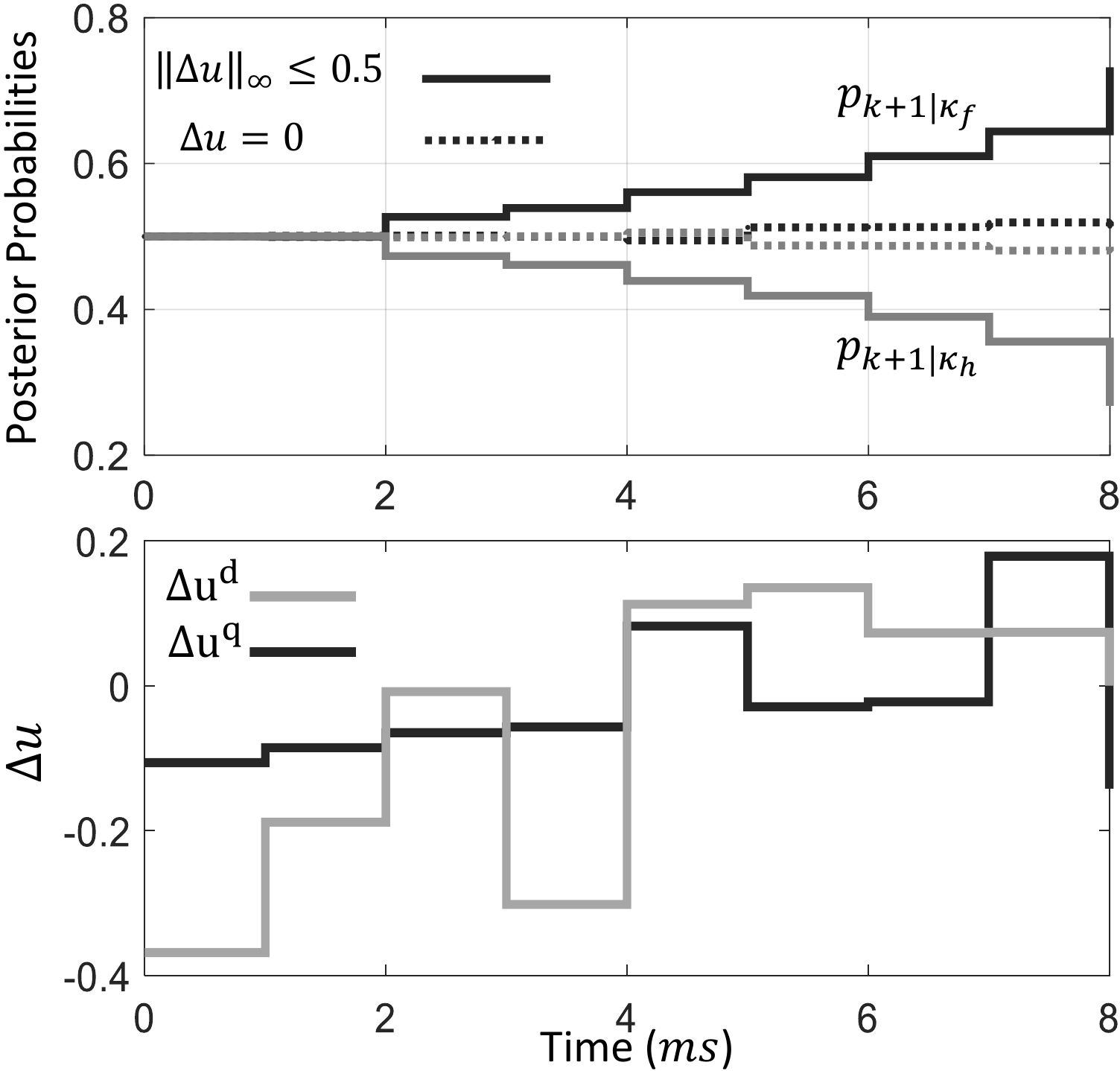}
\caption{(Top) Posterior probabilities with and without a perturbation. (Bottom) Perturbations over time in $d$ and $q$ frames.}
\label{fig:exasfvagsg3}
\end{figure}

In the following example, we see that without triggering the system with some additional perturbation, this method can fail to accurately detect faults. 
\begin{example}\label{exp:1}
Suppose that system \eqref{eqn:cna082} is faulty. For a  horizon of eight time steps, we calculate the the posterior probabilities \eqref{eqn:posterior} for the faulty and fault-free cases, starting from  $p_{0|f}=p_{0|h}=0.5$. In a second instance, we add a perturbation signal with a bounded norm, $\|\Delta u\|_{\infty}\leq 0.5$, to  help distinguish faulty from fault-free modes (the design of the input will be discussed in the next section). Fig.~\ref{fig:exasfvagsg3}  shows the evolution of $p_{k+1|i}$ over time. The passive method, i.e., with no perturbation, is far less sensitive to faults than the active method.
\end{example}

\subsection{Designing the Input Sequence}
As shown in Example \ref{exp:1},  the passive MMKF can lack adequate sensitivity on its own. We now design the perturbation, $\Delta u$, as shown in Fig.~\ref{fig:exwsghrle}, which makes the faulty and fault-free modes easily distinguishable from the system's output. More precisely, we consider the following problem.

\begin{problem}
Consider system \eqref{eqn:cna082} with $N>0$ as the detection horizon. Find an input sequence $\Delta u_k$, $k=0,1,2,...,N-1$, such that (i) $\|\Delta u_k\|_{\infty}\leq \gamma$, where  $\|\cdot\|_{\infty}$ is the infinity norm of a vector,  and $\gamma$ is a design parameter, and (ii) controller performance is not degraded, and (iii) the system's true mode is almost surely identified. The detection horizon $N$ and the bound on the input signal $\gamma$ are chosen by the designer.
\end{problem}

Note that the control input $u_k$ in \eqref{eqn:cna082} is augmented with the perturbation $\Delta u_k$. We  use the  objective function 
\begin{equation}\label{eqn:vn0a89}
J_d(\Delta u_{0:N-1})=\mathbb{E}\big[\sigma(\hat{\kappa})|\Delta u_{0:N-1}\big],
\end{equation}
where $\Delta u_{0:N-1}=[\Delta u_0~\cdots~\Delta u_{N-1}]\in\mathbb{R}^{N}$, $\hat{\kappa}$ is the identified mode of the system, and  $\sigma(\hat{\kappa})$ is zero when the identified mode is the actual
mode of the system, and is non-zero (e.g., equal to one) otherwise. In order to find the optimal input sequence for \eqref{eqn:vn0a89}, $J_d$ must be written in terms of $\Delta u_{0:N-1}$. We will show that \eqref{eqn:vn0a89} can be upper bounded by an explicit function of the input sequence $\Delta u_{0:N-1}$. To find that, we first need to calculate the average measurement $\bar{y}_{0:N|i}$ and its covariance matrix $\Sigma_{y|i}$ within the detection horizon. The mean output $\bar{y}_{0:N|i}$ is 
\begin{equation}
    \bar{y}_{0:N|i}=C_{i}\bar{x}_{0:N|i},
    \label{eqn:outpput}
\end{equation}
where $\bar{x}_{0:N|i}=[\bar{x}_{0|i}^T; \hspace{1mm} \bar{x}_{1|i}^T; \hspace{1mm}... \hspace{1mm} ;\bar{x}_{N|i}^T]^T$, and $ \bar{x}_{k|i}$ is the mean value of the state at time $k$,
\begin{equation}
    \bar{x}_{k|i}=A_{i}^k \bar{x}_{0|i}+ A_{i}^{k-1}B_{i}u_0+...+B_{i}u_{k-1}.
\end{equation}

The covariance matrix of the state at times $k$ and $l$ ($k\geq l$) is 
\begin{align}
\Sigma_{x,(k,l)|i}&=\mathbb{E}\left\{(x_{k|i}-\bar{x}_{k|i})(x_{k|i}-\bar{x}_{k|i})^T\right\}\nonumber\\
&=A_{i}^k\Sigma_{0}(A_{i}^l)^T+A_{i}^{k-1}\Sigma_{w}(A_{i}^{l-1})^T\nonumber\\
&+A_{i}^{k-2}\Sigma_{w}(A_{i}^{l-2})^T+...+A_{i}^{k-l}\Sigma_{w},
\end{align}
where $\Sigma_{0}$ and $\Sigma_{w}$ are covariance matrices of the initial state vector and the process noise, respectively. Due to the symmetry of covariance matrix, we have $\Sigma_{x,(k,l)|i}=\Sigma_{x,(l,k)|i}^T$.  The covariance matrix of the output is 
\begin{align}
\Sigma_{y,(k,l)|i}=
\begin{cases}
C_{i}\Sigma_{x,(k-1,l-1)|i}C_{i}^T \quad\quad\quad\hspace{1.5mm} k>l \\
C_{i}\Sigma_{x,(k-1,l-1)|i}C_{i}^T+\Sigma_{v} \quad k=l
\end{cases},
\label{eqn:22}
\end{align}
where $\Sigma_{v}$ is the covariance matrix of the measurement noise.

\begin{theorem}[\cite{Mehdi}]
The objective,  \eqref{eqn:vn0a89}, is bounded by
\begin{equation}\label{eqn:vn0a8sdg9}
J_d(\Delta u_{0:N-1})\leq \underbrace{\sqrt{p_{0|h}p_{0|f}}e^{-\phi}}_{\triangleq\hat{J}_d(\Delta u_{0:N-1})}
\end{equation}
where $p_{0|h}$ and $p_{0|f}$ and $\phi\geq0$ are prior probabilities, and $\phi\geq0$ is: 
\begin{align}\label{eq:phi}
\phi&=\frac{1}{4}\left(\bar{{y}}_{0:N|h}-\bar{y}_{0:N|f}\right)^T\left(\Sigma_{y|h}+\Sigma_{y|f}\right)^{-1}\nonumber\\
&\left(\bar{y}_{0:N|h}-\bar{y}_{0:N|f}\right)+\frac{1}{2}\ln\bigg(\frac{\det\left(\frac{\Sigma_{y|h}+\Sigma_{y|f}}{2}\right)}{\sqrt{\det\Sigma_{y|h}\det\Sigma_{y|f}}}\bigg).
\end{align}
\end{theorem}

In order to limit  the degradation of the control performance due to the addition of the perturbation $\Delta u_k$, we add a constraint on the norm of the input signal over the horizon. The final optimization takes the following form
\begin{equation}
    \begin{aligned}
        &\underset{\Delta u_{0:N-1}}{\text{min}}&& \hat{J}_d(\Delta u_{0:N-1})\\
        &\text{subject to} && \|\Delta u_{0:N-1}\|_{\infty}\leq \gamma.
    \end{aligned}
    \label{eqn:sidugb}
\end{equation}

\begin{remark}
On the one hand, since $p_{0|h}$ and $p_{0|f}$ are constant, minimizing $\hat{J}_d(\Delta u_{0:N-1})$ is equivalent to maximizing $\phi$. On the other hand, as seen in \eqref{eq:phi}, $\phi$ depends on the distance between expected outputs of fault-free and faulty systems. Therefore, minimizing $\hat{J}_d(\Delta u_{0:N-1})$ maximizes the distance between expected outputs of the fault-free and faulty systems. This means that the perturbation signal should be designed such that the distance between expected outputs of the fault-free and faulty systems is maximized. 
\end{remark}


\begin{remark}
In order to avoid degrading the system's performance in fault-free times, the perturbation $\Delta u$ is applied only when the reference inverter current $I_i^{\rm ref}$ approaches the upper limit of the FCL. Thus, the measure $m=|I_i^{\rm ref}-\zeta_U|$ is used as a passive indicator of fault detection. Note that relying on the difference between the reference current $I_i^{\rm ref}$ and the upper limit $\zeta_U$ may lead to a false alarm, as there are several scenarios in which $I_i^{\rm ref}$ can approach $\zeta_U$ without a fault in the system. For instance, an abrupt change in the load can lead to a large overshoot in the current reference generated by the voltage controller, and consequently saturate the FLC.
\end{remark}

\section{Simulations}
\label{sec:cm09n72}
We now evaluate the performance of the MMKF algorithm in simulation.  The parameters used for the voltage and current controllers are $k_{p}^I=170$, $k_{i}^I=100$
and $k_{p}^V=0.1$, $k_{i}^V=8$. The system parameters are $R=10\Omega$, $R_1=1.5m\Omega$, $L_1=300 mH$, and $V_{\rm DC}=150 V$. These physical parameters were taken from \cite{yazdani}.

\subsection{Detection-Performance Tradeoff}
\begin{figure}[t!]
\centering
\includegraphics[width=8.5cm]{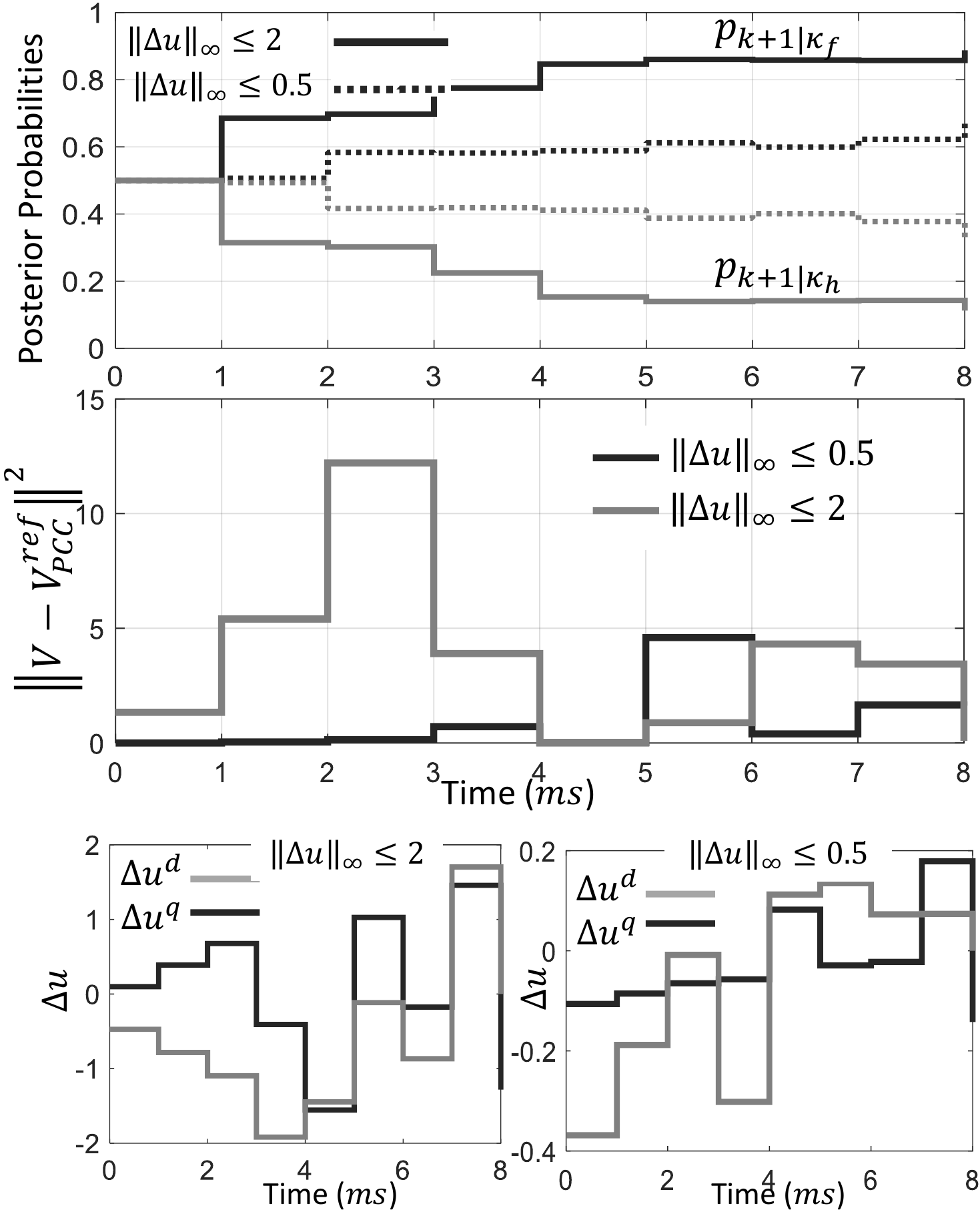}
\caption{The effect of the magnitude of the input signal $\Delta u$ on the detection time and tracking error for the faulty system}.
\label{fig:exnva98sg3}
\end{figure}

\begin{figure}[t!]
\centering
\includegraphics[width=8.5cm]{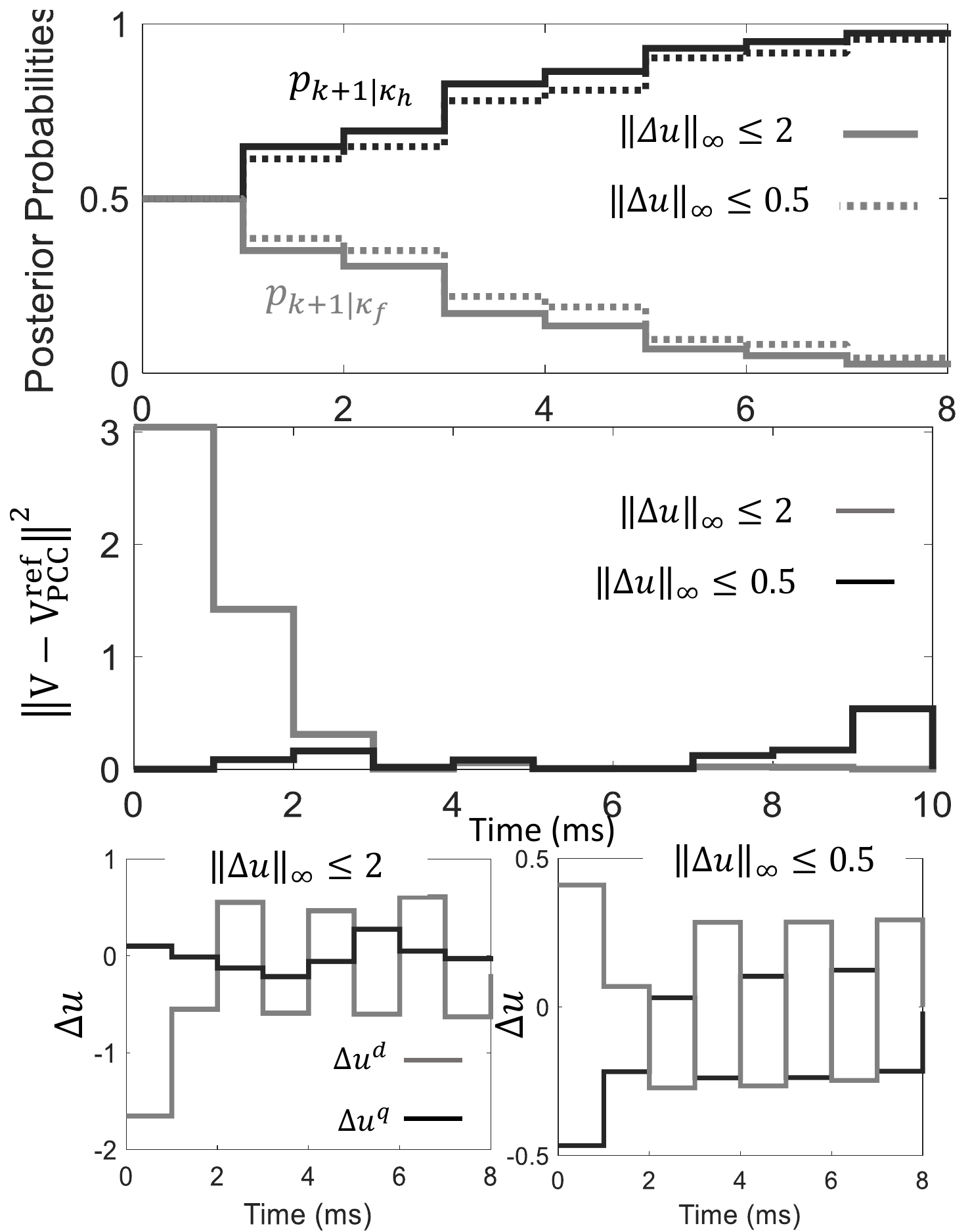}
\caption{The effect of the magnitude of the input signal $\Delta u$ on the detection time and tracking error for the fault-free system.}
\label{fig:FautFreeSim}
\end{figure}

As shown in Example \ref{exp:1}, perturbing the input improves fault detection. However, large perturbations may degrade performance. The design parameter $\gamma$ in \eqref{eqn:sidugb} plays an important role in balancing this tradeoff. Figs.~\ref{fig:exnva98sg3} and \ref{fig:FautFreeSim} show the effect of the perturbation $\Delta u_{0:N-1}$ on the deviation of the output voltage from its reference for the faulty and fault-free systems, respectively. As seen in these figures, faster and more accurate fault detection, quantified in terms of the separation of the posterior probabilities, comes at the price of performance degradation. We also observe that the separation of the posterior probabilities increases with time. More specifically, in both faulty and fault-free systems, the actual mode of the system can be detected within a single 60~Hz cycle.

\subsection{Effect of perturbation on the Output Voltage}
\begin{figure}[t!]
\centering
\includegraphics[scale=0.52]{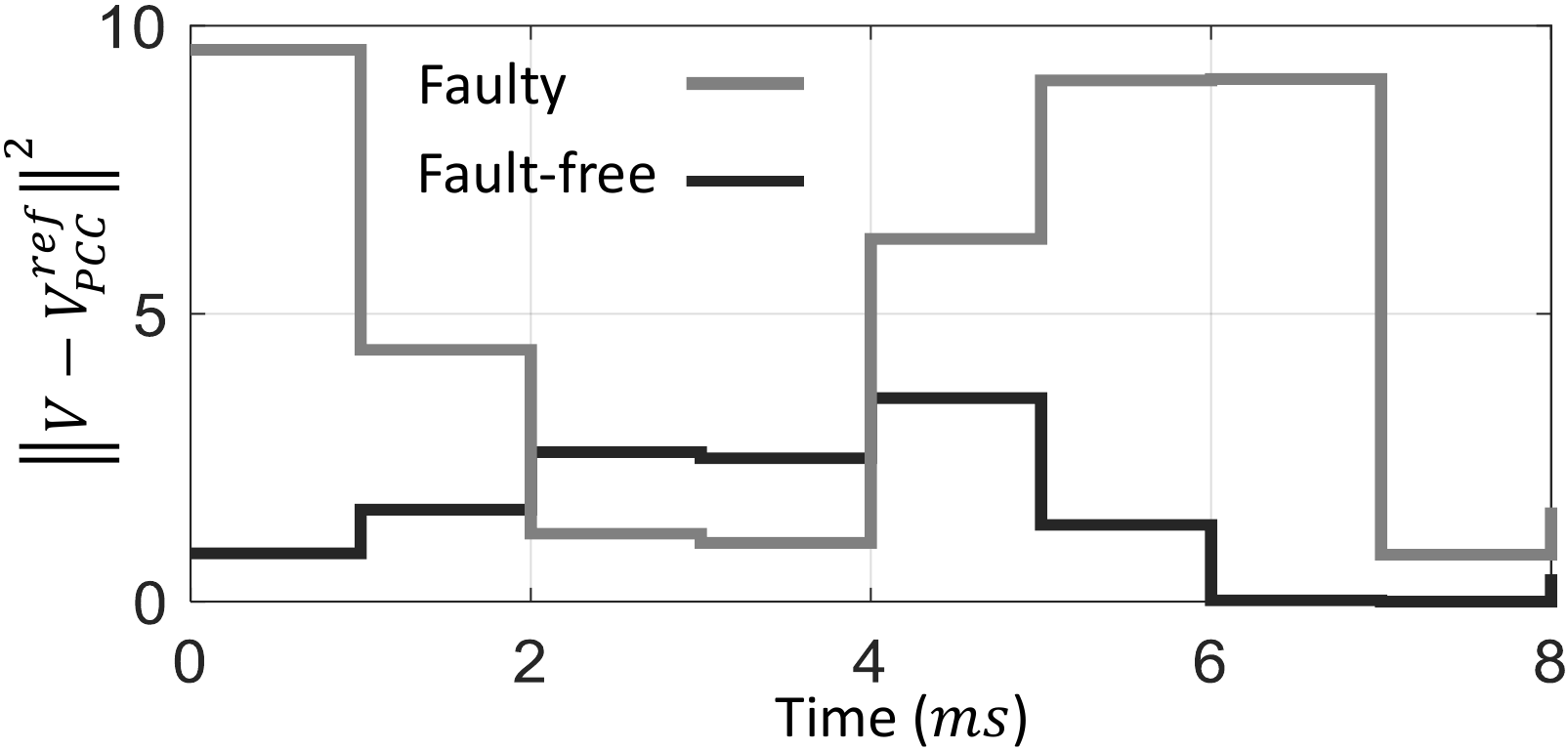}
\caption{The deviation from reference voltage for faulty and fault-free cases. }
\label{fig:exasafhorsg3}
\end{figure}

The influence of the perturbation on the output voltage can also be used as an indicator of faults. In the fault-free case, the controller consists of two cascaded PID blocks, while there is only one block in the faulty case. As a result, the perturbation will be more attenuated during fault-free than faulty operation. This is shown in Fig.~\ref{fig:exasafhorsg3}, in which a perturbation $\Delta u$ with $\|\Delta u_{0:N-1}\|_{\infty}\leq 1$ is applied to both the faulty and fault-free systems. The impact on the output voltage in faulty mode is significantly higher than in the fault-free case. 

\subsection{Comparison with Harmonic Injection}
\begin{figure}[t!]
\centering
\includegraphics[scale=0.55]{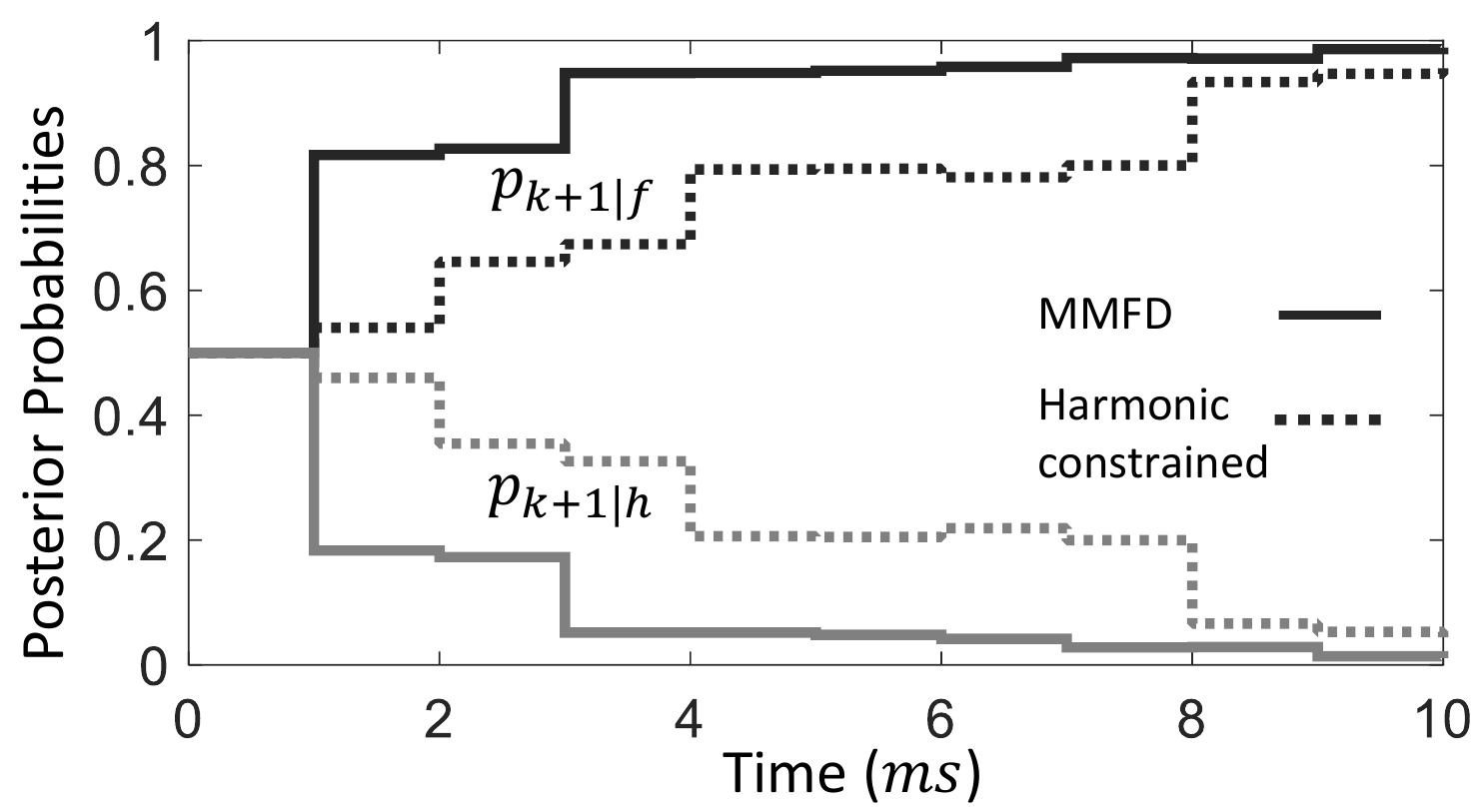}
\caption{The posterior probabilities for MMKF without the harmonics constraint (bold lines) and with the harmonic constraint (dashed lines). }
\label{fig:1698vdarsg3}
\end{figure}

We now compare our approach to \cite{khaled}, in which the perturbation is constructed from a few harmonics. We expect our perturbation to lead to better performance because it is optimized over a larger set of signals. 
 This is shown in Fig.~\ref{fig:1698vdarsg3}. The dashed lines correspond to the case where the input sequence $\Delta u_{0:N-1}$ is restricted to contain $3^{\rm rd}$, $5^{\rm th}$, and $7^{\rm th}$ harmonics of the fundamental frequency, 60~Hz. As expected, the probability of fault detection converges faster to 1 when the perturbation is not restricted to harmonics. 

\subsection{Robustness to Parameter Uncertainty}
An advantage of the MMKF is its robustness to model uncertainty. This is due to the fact that  uncertainties and disturbances are identically modeled in both the faulty and fault-free modes, and hence have a similar effect on the residual functions (and consequently on the posterior probabilities).
This feature is  beneficial, because identifying the exact parameters, and particularly the load $R$, is unrealistic in practice.

\begin{figure}[t!]
\centering
\includegraphics[scale=0.55]{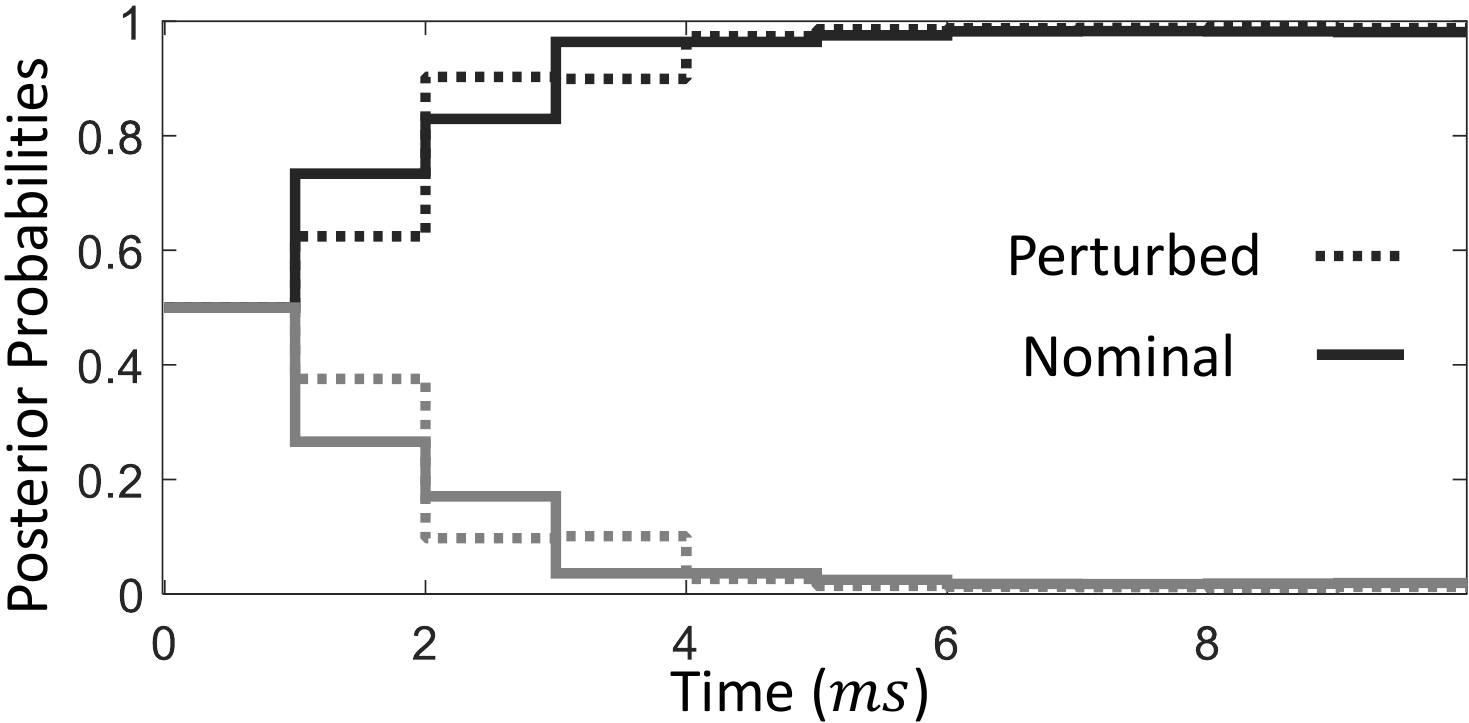}
\caption{The posterior probabilities for MMKF for nominal control parameters (bold line) and when they are perturbed by $10\%$ (dashed line). }
\label{fig:fmn0n1}
\end{figure}

As an example, we changed the parameters of the inner current controller by 10\%. As shown in Fig.~\ref{fig:fmn0n1}, this causes no significant change in the  posterior probabilities. 

\begin{figure}[t!]
\centering
\includegraphics[scale=0.58]{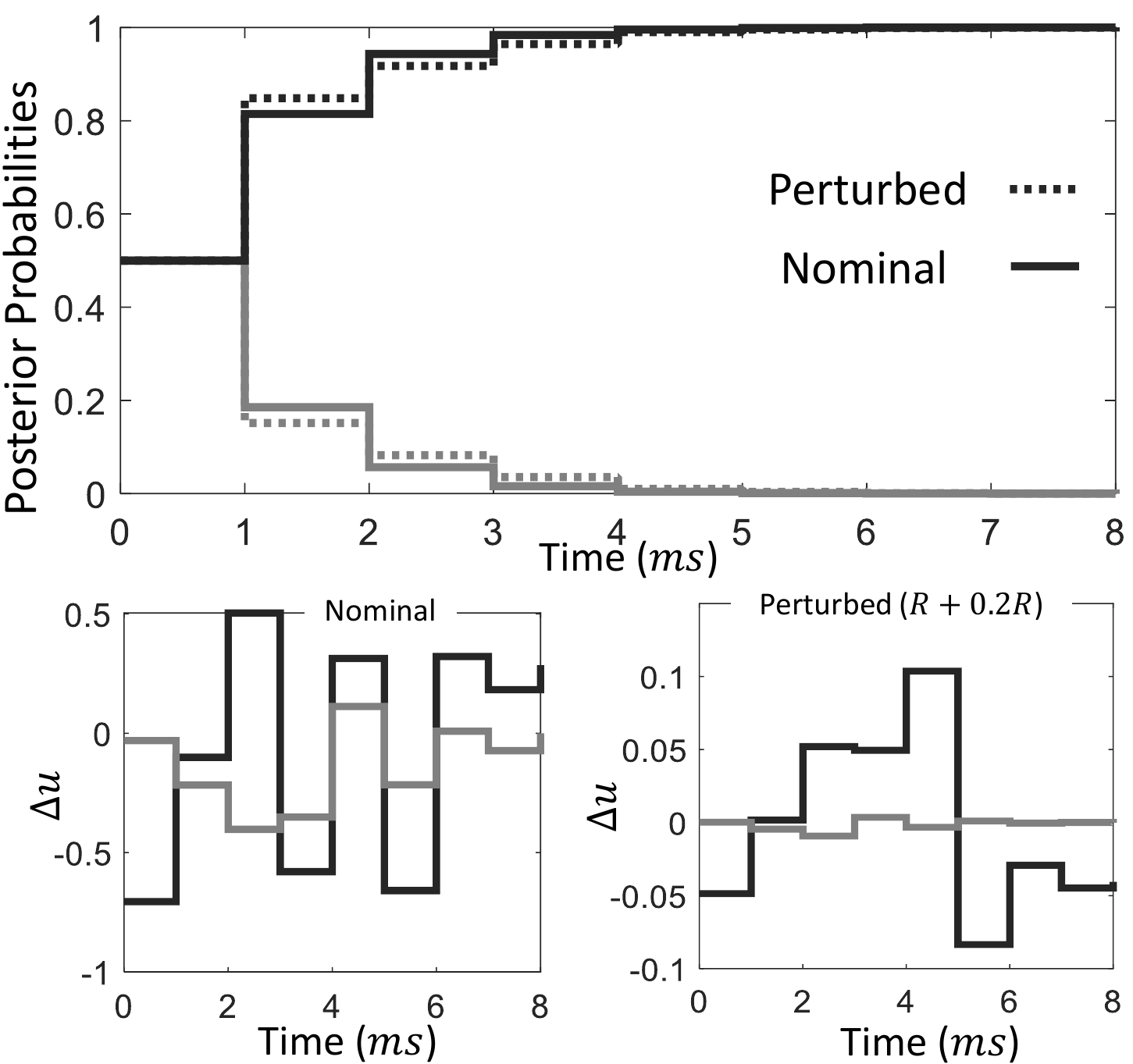}
\caption{(Top) The posterior probabilities for MMKF for a nominal load (bold line) and when it is perturbed by $20\%$ (dashed line). (Bottom) perturbation signal $\Delta u$ for nominal and perturbed load.  }
\label{fig:fmn0san1}
\end{figure}

To analyze the robustness of the proposed method against parameter uncertainty, we consider an identification error of $20\%$ in the load $R$ (i.e., $R\pm 0.2R$). Simulation results are shown in  Fig.~\ref{fig:fmn0san1}. Based on this figure, the algorithm can detect faults despite the existence of large uncertainties in the load. 
The perturbation inputs are shown in this figure. Based on this, when the load is increased ($R+ 0.2R$), the magnitude of the perturbations decreases. This is because when the load is larger, the difference between faulty (open loop, $R=0$) and fault-free cases becomes more apparent and less detection effort is needed.

\subsection{Sensitivity Analysis\textemdash Impact of Detection Horizon $N$}
The simulations are run on an Intel(R) Core(TM) i7-7500U CPU 2.70 GHz with 16.00 GB of RAM, and we use the YALMIP toolbox to solve the associated optimization problems. Setting $N=8$ (which is the length of the detection horizon in the results of the previous subsections), the mean computation time for 1000 runs is 0.612 ms. This is within the time available for real-time implementation, as the discretizing step is 1 ms.

\begin{figure}[t!]
\centering
\includegraphics[width=8.5cm]{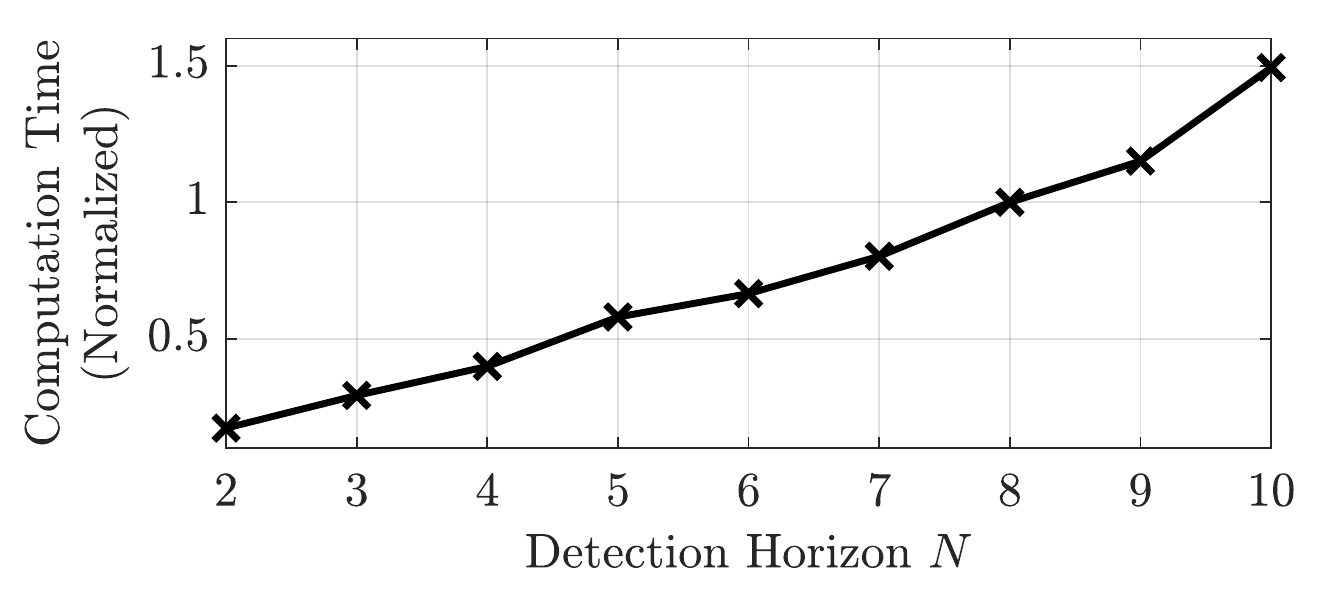}
\caption{Impact of detection horizon $N$ on the computation time, where the mean computation time for $N=8$ is used as normalizing coefficient.}
\label{fig:ComputationTime}
\end{figure}

Fig. \ref{fig:ComputationTime} shows how the detection window size impacts the computation time, where the mean computation time for $N=8$ (i.e., 0.612 ms) is used as a nominal value for normalization. This figure reports the mean computation time of 1000 experiments with random initial conditions for every $N$. Fig. \ref{fig:ComputationTime}  shows that as the length of the detection horizon increases, the optimization problem becomes larger, and the computation time increases.

\subsection{Sensitivity Analysis\textemdash Impact of Controller Parameters}
\begin{figure}[t!]
\centering
\includegraphics[width=8.5cm]{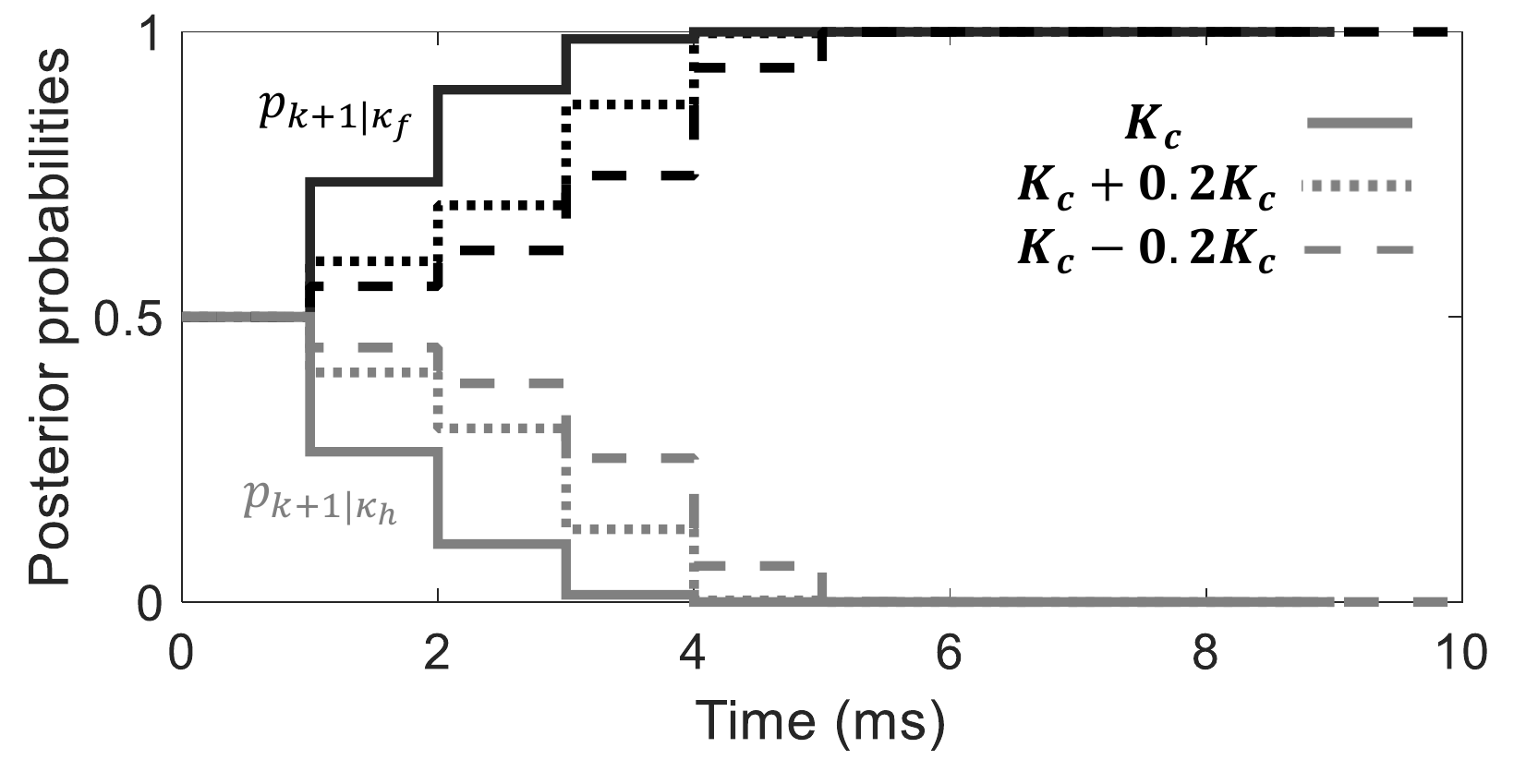}
\caption{Impact of PI controller parameters on detection performance.}
\label{fig:ControllerParameters}
\end{figure}

We follow the procedure in \cite{Hosseinzadeh2022} to maximize the distance between the faulty and fault-free system, which is assessed by means of the Vinnicombe gap (a.k.a. $\nu$-gap metric) \cite{Mahdianfar2011}. This procedure ensures that faulty and fault-free systems are far apart, which implies that the posterior probabilities converge fast. This is shown in Fig.~\ref{fig:ControllerParameters}, where $\boldsymbol K_c=[k_p^V~k_i^V~k_p^I~k_i^I]^\top$. As seen in this figure, even though we detect the actual mode of the system for all sets of controller parameters, the posterior probabilities for $\boldsymbol K_c$ converge faster than those for $\boldsymbol K_c\pm0.2 \boldsymbol K_c$.

\section{Summary and Conclusions}
 \label{sec:conclusion}
We have designed a new active scheme for detecting ground faults in inverter-based grids. The scheme consists of a multiple model Kalman filter together with an optimized perturbation, which enhances the sensitivity of the detector.  Simulations confirm that the detector is highly effective, that the perturbation signal significantly enhances detection, and that it outperforms harmonic injection-based techniques. Future work will focus on a broader range of faulty cases, e.g., unbalanced systems and single phase-to-ground fault.


\bibliographystyle{IEEEtran}
\bibliography{refs}
\end{document}